\documentclass[10pt]{article}
\usepackage{graphicx}
\usepackage{amssymb}
\usepackage{amsthm}
\usepackage{amsmath}
\usepackage{epsfig}
\def\ZZ{{\mathbf Z}}

\def\beqns{\begin{eqnarray*}}
\def\eeqns{\end{eqnarray*}}
\def\beqn{\begin{eqnarray}}
\def\eeqn{\end{eqnarray}}

\newcommand{\EQ}[2]{\begin{equation}{{#2}\label{#1}} \end{equation}}
                        \newcounter{theorem}

\newcommand{\newsection}[1]{{\setcounter{theorem}{0}%
                        \setcounter{equation}{0}}\section{#1}}
\numberwithin{equation}{section}
\numberwithin{theorem}{section}

\newenvironment{thm}[2]{\begin{sloppypar}\refstepcounter{theorem}%
                        {\bf #1 \thetheorem.}\label{#2}\em{}}%
                        {\end{sloppypar}}
                        \newcommand{\theo}[3]{\begin{thm}{#1}{#2}
#3\end{thm}}

                        \newcommand{\N}{{\rm I}\!{\rm N}}

                        \newcommand{\R}{{\rm I}\!{\rm R}}
\begin{document}

\title{Convergence of sparse grid Gaussian convolution approximation for multi-dimensional periodic functions}
\author{Simon Hubbert\footnote{Department of Economics, Mathematics and Statistics, Birkbeck, University of London, WC1H 7HX, UK. {\tt s.hubbert@bbk.ac.uk}}, 
Janin J\"ager\footnote{Lehrstuhl Numerische Mathematik, Justus-Liebig University, 35392 Giessen, Germany. {\tt janin.jaeger@math.uni-giessen.de}} and Jeremy Levesley\footnote{Department of Mathematics, University of Leicester, LE1 7RH, UK. {\tt jl1@le.ac.uk}} }

\maketitle

\begin{abstract}
We consider the problem of approximating $[0,1]^{d}$-periodic functions by convolution with a scaled Gaussian kernel. We start by establishing convergence rates to functions from periodic Sobolev spaces and we  show that the saturation rate is $O(h^{2}),$ where $h$ is the scale of the Gaussian kernel. Taken from a discrete point of view, this result can be interpreted as the accuracy that can be achieved on the uniform grid with spacing $h.$ In the discrete setting, the curse of dimensionality would place severe restrictions on the computation of the approximation. For instance, a spacing of $2^{-n}$ would provide an approximation converging at a rate of $O(2^{-2n})$ but would require  $(2^{n}+1)^{d}$ grid points. To overcome this we introduce a sparse grid version of Gaussian convolution approximation, where substantially fewer grid points are required, and
 show that the sparse grid version delivers a saturation rate of $O(n^{d-1}2^{-2n}).$ This rate is in line with what one would expect in the sparse grid setting (where the full grid error only deteriorates by a factor of order $n^{d-1}$) however the analysis that leads to the result is novel in that it draws on results from the theory of special functions and key observations regarding the form of certain weighted geometric sums.

\end{abstract}

\section{Introduction}

Many methods that are designed to deliver approximations are based on the convolution of a kernel function with the function being approximated. The general approach involves selecting a suitable integrable function $K:\R^{d}\to \R$ ({\it{the convolution kernel}}) satisfying
\[
\int_{\R^{d}}K({\bf{x}})d{\bf{x}}=\int_{\R^{d}}K(x_{1},\ldots,x_{d})dx_{1} \cdots dx_{d}=1.\]
 A scaling vector ${\bf{h}}=(h_{1},\ldots,h_{d})^{T}\in \R^{d},$ with $h_{i}>0 $ ($1\le i\le d$) is then used to define a parameterised family of convolution kernels by
\[K_{{\bf{h}}}({\bf{x}})=\frac{1}{h_{1}\cdots h_{d}}K\left(\frac{x_{1}}{h_{1}},\ldots,\frac{x_{d}}{h_{d}}\right).
\]
The convolution approximation to  a function $f \in L_{1}(\R^{d})$  is defined by
\EQ{aniso}{
\mathcal{C}_{\textbf{h}}(f)(\textbf{x})=(f \ast K_{{\bf{h}}})({\bf{x}})=\int_{\R^{d}}f({\bf{y}})K_{{\bf{h}}}({\bf{x}}-{\bf{y}})d{\bf{y}}.
}

The  convolution kernel described above is anisotropic as each direction is allowed to be scaled by its own factor. This is practically useful because a typical data sample of a function will show variety along different directions and so a well designed anisotropic scaling can efficiently capture these features. However,  from a theoretical perspective, most convergence results relate to the isotropic scale where each direction is scaled by the same factor $h>0.$    In this case the scaled kernel is $K_{h}({\bf{x}})=h^{-d}K\left(\frac{{\bf{x}}}{h}\right)$ and the corresponding convolution approximation
\EQ{iso}{
\begin{aligned}
{\cal{C}}_{h}(f)({\bf{x}})=\int_{\R^{d}}f({\bf{y}})K_{h}({\bf{x}}-{\bf{y}})d{\bf{y}}=\frac{1}{h^{d}}\int_{\R^{d}}f({\bf{y}})K\left(\frac{{\bf{x}}-{\bf{y}}}{h}\right)d{\bf{y}},
\end{aligned}
}
can  be shown to converge to $f$ as $h \to 0$, the convergence being uniform on compact sets, \cite{cheneylight} chapter 20, theorem 2. The rate of convergence depends upon the smoothness of $f$ and the polynomial reproduction properties of the underlying kernel.
The  convolution approximation can be viewed as the continuous counterpart of quasi-interpolation; a discrete method which generates an approximation over the whole of $\R^{d}$ by linearly combining the values of $f$ sampled at the scaled integer lattice $h{{\ZZ}}^{d}$ together with the appropriately shifted and scaled kernel function. The classical construction, as for example described in \cite{Buhmann}, takes the form
\EQ{quasi}{
Q_{h}(f)({\bf{x}}) = \sum_{{\bf{k}}\in {{\ZZ}}^{d}}f({\bf{k}}h)K\left(\frac{{\bf{x}}}{h}-{\bf{k}}\right),\quad {\bf{x}}\in \R^{d}, \quad h>0.
}

Following \cite{mazya} the connection between continuous convolution and discrete quasi-interpolation can be seen if we write
\[
\begin{aligned}
{\cal{C}}_{h}(f)({\bf{x}})=\frac{1}{h^{d}}\sum_{{\bf{k}}\in {\ZZ}^{d}}\int_{h\cdot({\bf{k}}+[-\frac{1}{2},\frac{1}{2}]^{d})}f({\bf{y}})K\left(\frac{{\bf{x}}-{\bf{y}}}{h}\right)d{\bf{y}}.
\end{aligned}
\]
The integrals above are taken over appropriately shifted and scaled versions of the cube $[-\frac{1}{2},\frac{1}{2}]^{d}.$ If we approximate each integrand by its value at the midpoint of the cube we get
\[
  \int_{h\cdot({\bf{k}}+[-\frac{1}{2},\frac{1}{2}]^{d})}f({\bf{y}})K\left(\frac{{\bf{x}}-{\bf{y}}}{h}\right)d{\bf{y}}\approx h^{d}f({\bf{k}}h)K\left(\frac{{\bf{x}}}{h}-{\bf{k}}\right),
\]
and so we have that ${\cal{C}}_{h}(f)({\bf{x}}) \approx Q_{h}(f)({\bf{x}}).$ Quasi-interpolation using Gaussians in one dimension was described in \cite{JerHub}.

In this paper we will examine the approximation of $[0,1]^{d}-$periodic functions by convolution with the multi-dimensional Gaussian kernel. Given the close connection of continuous convolution to quasi-interpolation  the results we establish in the continuous setting will serve as a baseline for what should be expected in the discrete case. 

We begin in Section 2 by  deriving the formula for the Fourier expansion of the pointwise error using the anisotropic scaling of the Gaussian; this result allows us to deduce that convolution approximation  is only able to reproduce the constant function. We then analyse the isotropic case in some detail. In this setting we demonstrate that the convergence has a saturation rate of ${\cal{O}}(h^{2})$. \\

In Section 3 we consider the practical issues of employing the discrete (quasi-interpolation) analogue of continuous convolution in high dimensions. Such a recasting involves constructing a full grid in $[0,1]^{d}$ with an isotropic spacing of $h=\frac{1}{2^{n}},$ where $n$ is a positive integer. In this setting,  the convolution approximation will converge to $f$ at a rate of $\frac{1}{2^{2n}},$ provided $f$ is sufficiently smooth. However, in the discrete setting we are restricted by the curse of dimensionality since the construction of the quasi-interpolant would require $(2^{n}+1)^{d}$ evaluations and this is prohibitively large as $n$ grows. In order to overcome this we consider replacing the full-grid approximation with a sparse grid version which is built from a certain linear combination of smaller full grid approximations. Numerical experiments on closely connected methods have been published in \cite{UstaLev}. To analyse this theoretically we mimic the approach of Section 2, i.e., we first derive the formula for the Fourier expansion of the pointwise error using the sparse grid convolution approximation. We then investigate  the Fourier coefficients of the error expansion and we state the main theorem of the paper, concerning the decay rate of the coefficients. We then establish that, provided $f$ is sufficiently smooth,  the sparse grid convolution approximation will converge to $f$ at a rate of $\frac{n^{d-1}}{2^{2n}}.$  Section 4 is devoted to the proof of the aforementioned main theorem of the paper.

%
%

\section{Gaussian convolution approximation}

Our choice of convolution kernel is the multi-dimensional Gaussian
\[
\begin{aligned}
\Psi({\bf{x}})=\frac{1}{(2\pi)^{\frac{d}{2}}}\exp\left(-\frac{1}{2}{\bf{x}}^{T}{\bf{x}}\right)
=\frac{1}{(2\pi)^{\frac{d}{2}}}\exp\left(-\frac{1}{2}\left(\sum_{i=1}^{d}x_{i}^{2}\right)\right)=
\prod_{i=1}^{d}\psi(x_{i}),
\end{aligned}
\]
where $\psi:\R\to\R$ is the univariate Gaussian $\psi(x) = (2\pi)^{-\frac{1}{2}}\exp\left(-\frac{1}{2}x^{2}\right).$ Fourier theory will play an important role in our analysis and we recall that if we let $e_{x}(z)= \exp(2\pi i x\cdot z)$ then the univariate Fourier transform of $\psi$ is
\[
\widehat{\psi}(z):=\int_{-\infty}^{\infty}\psi(x)e_{-x}(z)dx=\exp(-2\pi^{2}z^{2}).
\]

Our general aim is  to approximate a $[0,1]^{d}$-periodic function
\[
f({\bf{x}})=\sum_{{\bf{k}}\in {\ZZ}^{d}}{\widehat{f}({\bf{k}})}e_{{\bf{k}}
}({\bf{x}})\quad
{\rm{where}}
\quad  e_{{\bf{y}}}({\bf{x}}):=\exp(2\pi i {\bf{y}}^{T}{\bf{x}})=\prod_{i=1}^{d}e_{y_{i}}(x_{i}),
\]
by the continuous multi-variable convolution
\[
\begin{aligned}
&\mathcal{C}_{{\bf{h}}}(f)
=\int_{\R^{d}} f(\textbf{z})\Psi_{\textbf{h}}(\textbf{x}-\textbf{z})d\textbf{z}=\sum_{\textbf{j}\in {\ZZ}^{d}}\int_{[0,1]^{d}}f(\textbf{z})\Psi_{\textbf{h}}(\textbf{x}-\textbf{z}-\textbf{j})
d\textbf{z}=\int_{[0,1]^{d}}f(\textbf{z})\Phi_{\textbf{h}}(\textbf{x}-\textbf{z})d\textbf{z},
\end{aligned}
\]
where
\[
\Phi_{\textbf{h}}(\textbf{x})=\sum_{\textbf{j}\in {\ZZ}^{d}}\Psi_{\textbf{h}}(\textbf{x}-\textbf{j}).
\]

Now, $\Phi_{\textbf{h}}(\textbf{x})$ is $[0,1]^{d}-$periodic and so has a multi-dimensional Fourier series
\[
\Phi_{\textbf{h}}(\textbf{x})=\sum_{\textbf{k}\in {\ZZ}^{d}}\widehat{\Phi_{\textbf{h}}}(\textbf{k})e_{\textbf{k}}(\textbf{x}),
\]
where
\[
\begin{aligned}
\widehat{\Phi_{\textbf{h}}}(\textbf{k})&=\int_{[0,1]^{d}}\Phi_{\textbf{h}}(\textbf{x})
e_{-\textbf{k}}(\textbf{x})d\textbf{x}\\
&=\int_{[0,1]^{d}}\left(\sum_{\textbf{j}\in {\ZZ}^{d}}\Psi_{\textbf{h}}(\textbf{x}-\textbf{j})\right)
e_{-\textbf{k}}(\textbf{x})d\textbf{x}\\
&=\sum_{\textbf{j}\in {\ZZ}^{d}}\int_{[0,1]^{d}}
\Psi_{\textbf{h}}(\textbf{x}-\textbf{j})e_{-\textbf{k}}(\textbf{x})d\textbf{x}\\
&=\int_{\R^{d}}\Psi_{\textbf{h}}(\textbf{x})e_{-\textbf{k}}(\textbf{x})d\textbf{x}
=\prod_{i=1}^{d}\int_{-\infty}^{\infty}\psi_{h_{i}}(x_{i})e_{-k_{i}}(x_{i})dx_{i}=\prod_{i=1}^{d}\widehat{\psi}(h_{i}k_{i}).
\end{aligned}
\]
Applying the $d$-dimensional convolution formula for $[0,1]^{d}-$periodic functions we have:
\EQ{genconapp}{
\begin{aligned}
\mathcal{C}_{{\bf{h}}}(f)(\textbf{x})=\sum_{\textbf{k}\in {\ZZ}^{d}}\widehat{f}(\textbf{k})\widehat{\Phi_{\textbf{h}}}(\textbf{k})e_{\textbf{k}}(\textbf{x})=\sum_{\textbf{k}\in {\ZZ}^{d}}\widehat{f}(\textbf{k})\left(\prod_{i=1}^{d}\widehat{\psi}(h_{i}k_{i})\right)e_{\textbf{k}}(\textbf{x}).
\end{aligned}
}
Thus the error in the convolution approximation is
\EQ{errorbasic}{
\begin{aligned}
E_{\textbf{h}}(f)({\bf{x}})=
f({\bf{x}})-\mathcal{C}_{{\bf{h}}}(f)(\textbf{x}) = \sum_{{\bf{k}}\in {\ZZ}^{d} }\widehat{f}(\textbf{k})\left(1-\left(\prod_{i=1}^{d}\widehat{\psi}(h_{i}k_{i})\right)\right)
e_{\textbf{k}}(\textbf{x}).
\end{aligned}}
We note that, since  $\widehat{\psi}(0)=1,$ the above error representation immediately shows that the convolution reproduces the constant but not any other trigonometric polynomial.

\subsection{Convergence with Isotropic scaling}

If we consider the isotropic case where the same scale factor $h$ is applied to all coordinate directions then (\ref{errorbasic}) can be written as
\EQ{basiciso}
{
E_{h}(f)({\bf{x}})= \sum_{{\bf{k}}\in {\ZZ}^{d} \setminus \{ {\bf{0}}\}}\widehat{f}(\textbf{k})\left(1-e^{-2\pi^{2}h^{2}\|{\bf{k}}\|^{2}}\right)e_{\textbf{k}}(\textbf{x}).}

The functions we wish to approximate are  taken from a  periodic Sobolev space
$$
\mathcal{N}_\beta = \left\{ f=\sum_{{\bf{k}}\in {\ZZ}^{d}} {\widehat{f}({\bf{k}})} e_{\bf{k}} : \| f \|_{\beta} = \left ( \sum_{\textbf{k}\in {\ZZ}^{d}\setminus\{{\textbf{0}}\}} \|\textbf{k}\|^{2\beta}|{\widehat{f}({\bf{k}})}|^{2} \right )^{1/2} < \infty \right\}.
$$
The Sobolev embedding theorem \cite{adams} ensures that if $\beta>\frac{d}{2}$ then all functions in $\mathcal{N}_\beta$ will be continuous. The following result gives error bounds for Gaussian convolution approximation of such functions.

\theo{Proposition}{prop1}{
Let $f \in \mathcal{N}_\beta,$ where $\beta >\frac{d}{2}.$ Then
\[
\| E_h f \|_{\infty}\le \|f\|_{\beta} \cdot \begin{cases} C_{1}h^2 \,\, & \,\,
\textrm{for} \,\,\,\, \beta>\frac{d}{2}+2; \\ h^2\left(C_{2}\sqrt{\ln\left(\frac{1}{h}\right)}+C_{3}\right) \,\, & \,\,
\textrm{for} \,\,\,\, \beta=\frac{d}{2}+2; \\
C_{4}h^{\beta-\frac{d}{2}} \,\, & \,\,
\textrm{for} \,\,\,\, \frac{d}{2}< \beta < \frac{d}{2}+2,
\end{cases}
\]
where $C_{i}$ $i=1,2,3,4$, are positive constants independent of $h.$
}
\begin{proof}
Using   (\ref{basiciso}) together with the elementary bound $ 1- e^{-x}<x$ (for $x>0),$ we can deduce that

\EQ{base-error}{
\| E_h f \|_{\infty
} \le 2\pi^{2}h^{2}\sum_{{\bf{k}}\in {\ZZ}^{d} \setminus \{ {\bf{0}}\}}\|{\bf{k}}\|^{2}|{\widehat{f}({\bf{k}})}|.
}

Suppose that  $\beta = \frac{d}{2}+2+\alpha,$ where $\alpha>0,$ then an application of
 the Cauchy Schwarz inequality yields
\[
\begin{aligned}
\| E_h f \|_{\infty}&\le 2\pi^{2}h^{2}\sum_{{\bf{k}}\in {\ZZ}^{d} \setminus \{ {\bf{0}}\}}\frac{1}{\|{\bf{k}}\|^{\frac{d}{2}+\alpha}}\cdot \|{\bf{k}}\|^{\frac{d}{2}+\alpha+2}|{\widehat{f}({\bf{k}})}|\\
&\le 2\pi^{2}h^{2}\left(\sum_{{\bf{k}}\in {\ZZ}^{d} \setminus \{ {\bf{0}}\}}
\frac{1}{\|{\bf{k}}\|^{d+2\alpha}}\right)^{\frac{1}{2}}\cdot \left(\sum_{{\bf{k}}\in {\ZZ}^{d} \setminus \{ {\bf{0}}\}}\left(\|{\bf{k}}\|^{2}\right)^{ \frac{d}{2}+2+\alpha}|{\widehat{f}({\bf{k}})}|^{2}\right)^{\frac{1}{2}}\\
&\le 2\pi^{2}C h^{2}\|f\|_{\beta}.
\end{aligned}
\]
Now assume that $\beta = \frac{d}{2}+2-\alpha$ where $0 < \alpha < 2.$ In the following development we will  work with a partition of the punctured integer lattice
\[
{\ZZ}^{d} \setminus \{ {\bf{0}}\}= \underbrace{\left\{{\bf{k}}\in {\ZZ}^{d} \setminus \{ {\bf{0}}\}: 1\le \|{\bf{k}}\|< \frac{1}{h}\right\}}_{=\Theta_{h}} \cup
\underbrace{\left\{{\bf{k}}\in {\ZZ}^{d} :\|{\bf{k}}\| \ge \frac{1}{h}\right\}}_{=\Gamma_{h}}.\]

Using this we bound the error in two parts as follows

\[
\| E_h f \|_{\infty} \le \sum_{\mathbf{k}\in \Theta_{h}}|{\widehat{f}({\bf{k}})}|\cdot|1-e^{-2\pi^{2}h^{2}\|\mathbf{k}\|^{2}}|
+\sum_{\mathbf{k}\in \Gamma_{h}}|{\widehat{f}({\bf{k}})}|\cdot|1-e^{-2\pi^{2}h^{2}\|\mathbf{k}\|^{2}}|.
\]
For the sum over $\Theta_{h}$ we again employ $ 1- e^{-x}<x$  and, bounding as before, we conclude that
\[
\begin{aligned}
&\sum_{\mathbf{k}\in \Theta_{h}}|{\widehat{f}({\bf{k}})}|\cdot|1-e^{-2\pi^{2}h^{2}\|\mathbf{k}\|^{2}}|
\le 2\pi^{2}h^{2}\sum_{\mathbf{k}\in \Theta_{h}}|{\widehat{f}({\bf{k}})}|\cdot\|\mathbf{k}\|^{2}\\
&=2\pi^{2}h^{2}\sum_{\mathbf{k}\in \Theta_{h}}\|\mathbf{k}\|^{\alpha-\frac{d}{2}}\|\mathbf{k}\|^{\frac{d}{2}+2-\alpha}|{\widehat{f}({\bf{k}})}|\\
&\le 2\pi^{2}h^{2} \left(\sum_{\mathbf{k}\in \Theta_{h}}\|\mathbf{k}\|^{2\alpha-d}\right)^{\frac{1}{2}}
\left(\sum_{\mathbf{k}\in \Theta_{h}}\left(\|\mathbf{k}\|^{2}\right)^{\frac{d}{2}+2-\alpha}|{\widehat{f}({\bf{k}})}|^{2}
\right)^{\frac{1}{2}}\\
&\le 2\pi^{2}h^{2} \left(\sum_{\mathbf{k}\in \Theta_{h}}\left(\frac{1}{h}\right)^{2\alpha-d}\right)^{\frac{1}{2}}\|f\|_{\frac{d}{2}+2-\alpha}
\\
&\le C\cdot 2\pi^{2}h^{2}\Bigl[\left(\frac{1}{h}\right)^{d}\left(\frac{1}{h}\right)^{2\alpha-d}\Bigr]^{\frac{1}{2}}\|f\|_{\frac{d}{2}+2-\alpha}\\
&\le C\cdot 2\pi^{2}h^{2-\alpha}\|f\|_{\frac{d}{2}+2-\alpha}=C\cdot 2\pi^{2}h^{\beta-\frac{d}{2}}\|f\|_{\beta}.
\end{aligned}
\]

We note that for the case where $\alpha = 0,$ corresponding to $\beta = \frac{d}{2}+2,$ the above development can be traced to the third line to yield
\[
\sum_{\mathbf{k}\in \Theta_{h}}|{\widehat{f}({\bf{k}})}|\cdot|1-e^{-2\pi^{2}h^{2}\|\mathbf{k}\|^{2}}|
\le 2\pi^{2}h^{2}\left(\sum_{\mathbf{k}\in \Theta_{h}}\frac{1}{\|\mathbf{k}\|^{d}}\right)^{\frac{1}{2}}\|f\|_{\beta}.
\]

Applying the integral  test, with a change to polar coordinates, we have
\[
  \sum_{\mathbf{k}\in \Theta_{h}}\frac{1}{\|\mathbf{k}\|^{d}}\le \int_{1\le \|{\bf{x}}\|\le \frac{1}{h}} \frac{d{\bf{x}}}{\|{\bf{x}}\|^{d}}=C_{d}\int_{1}^{\frac{1}{h}}\frac{dr}{r}=C_
  {d}\ln\left(\frac{1}{h}\right).
  \]

In summary, for this part of the sum we can conclude that
\EQ{partone}{
\sum_{\mathbf{k}\in \Theta_{h}}|{\widehat{f}({\bf{k}})}|\cdot|1-e^{-2\pi^{2}h^{2}\|\mathbf{k}\|^{2}}| \le  \|f\|_{\beta}
\begin{cases} Ch^2\sqrt{\ln\left(\frac{1}{h}\right)} \,\, & \,\,
\textrm{for} \,\,\,\, \beta=\frac{d}{2}+2; \\
Ch^{\beta-\frac{d}{2}} \,\, & \,\,
\textrm{for} \,\,\,\, \frac{d}{2}< \beta < \frac{d}{2}+2.
\end{cases}
}

For the sum over $\Gamma_{h}$ we employ $ 1- e^{-x}< 1$ and develop the bound as follows:

\[
\begin{aligned}
&
\sum_{\mathbf{k}\in \Gamma_{h}}|{\widehat{f}({\bf{k}})}|\cdot|1-e^{-2\pi^{2}h^{2}\|\mathbf{k}\|^{2}}|
\le \sum_{\mathbf{k}\in \Gamma_{h}}|{\widehat{f}({\bf{k}})}|=\sum_{\mathbf{k}\in \Gamma_{h}}|\|\mathbf{k}\|^{\alpha-2-\frac{d}{2}}\cdot \|\mathbf{k}\|^{\frac{d}{2}+2-\alpha}|{\widehat{f}({\bf{k}})}|\\
&\le \left(\sum_{\mathbf{k}\in \Gamma_{h}}\|\mathbf{k}\|^{2\alpha-4-d}\right)^{\frac{1}{2}}\left(\sum_{\mathbf{k}\in \Gamma_{h}}\left(\|\mathbf{k}\|^{2}\right)^{\frac{d}{2}+2-\alpha}|{\widehat{f}({\bf{k}})}|^{2}
\right)^{\frac{1}{2}}\\
&\le \left(\sum_{\mathbf{k}\in \Gamma_{h}}\|\mathbf{k}\|^{2\alpha-4-d}\right)^{\frac{1}{2}}\|f\|_{\frac{d}{2}+2-\alpha}=\left(\sum_{\mathbf{k}\in \Gamma_{h}}\|\mathbf{k}\|^{-2\beta}\right)^{\frac{1}{2}}\|f\|_{\beta}.
\end{aligned}
\]

Keeping in mind that $\frac{d}{2}<\beta  \le
 \frac{d}{2}+2,$   the integral comparison test  yields
\[
 \sum_{\mathbf{k}\in \Gamma_{h}}\|\mathbf{k}\|^{-2\beta}\le \int_{\|{\bf{x}}\|\ge \frac{1}{h}}\|{\bf{x}}\|^{-2\beta}d{\bf{x}}=C_{d}\int_{\frac{1}{h}}^{\infty}r^{d-1-2\beta}dr\le Ch^{2\beta-d},
 \]
 and, for this parameter range, we   can deduce that
 \EQ{parttwo}{
 \sum_{\mathbf{k}\in \Gamma_{h}}|{\widehat{f}({\bf{k}})}|\cdot|1-e^{-2\pi^{2}h^{2}\|\mathbf{k}\|^{2}}|\le
 Ch^{\beta-\frac{d}{2}},\,\,\,{\rm{for}}\,\,\, \frac{d}{2}<\beta  \le
 \frac{d}{2}+2.
 }
Combining this with (\ref{partone}) provides the bounds stated in the proposition.

\end{proof}

\section{Gaussian convolution approximation on sparse grids}

The convergence results of the previous section are of theoretical interest, however, from a practical perspective, the implementation of the discrete (quasi-interpolation) analogue in high dimensions is restricted by the curse of dimensionality. A direct recasting of the continuous case to discrete setting would require that we sample values on a full grid in $[0,1]^{d},$ thus for $h=1/2^{n}$ this would amount to  $(2^{n}+1)^{d}$ evaluations.   One remedy that can be used to alleviate the curse of dimensionality, at least for moderately high dimensions,  is to approximate on a carefully chosen subset of the full grid, where substantially fewer points are needed to achieve an acceptable level of accuracy. To describe our approach we let
 $\boldsymbol\ell= (\ell_{1},\ldots,\ell_{d})$ denote a general multi-index where $\ell_{i} \ge 1$ for $i=1,\ldots,d$ then we define $\mathcal{X}_{\boldsymbol\ell}$ to be the anisotropic (directionally uniform) grid in $[0,1]^{d}$ where $h_{i}=1/2^{\ell_{i}}$ denotes the spacing in the $i^{th}$ coordinate direction. The number of nodes in $\mathcal{X}_{\boldsymbol\ell}$  is then given by
\[
|\mathcal{X}_{\boldsymbol\ell}| = \prod_{i=1}^{d}(2^{\ell_{i}}+1).
\]

We let $\mathcal{X}_{n,d}=\mathcal{X}_{2^{-n}{\bf{1}}}$ denote the full isotropic grid with a uniform spacing of $h = \frac{1}{2^{n}}.$ As a starting point we can appeal to Proposition \ref{prop1} to conclude that the approximation error for the continuous convolution approximation to any given $f \in \mathcal{N}_\beta$ $(\beta >\frac{d}{2}+2)$ on the full grid $\mathcal{X}_{n,d}$ satisfies

\EQ{full}{
\|E_{\frac{1}{2^{n}}}f\|_{\infty} =\|f-\mathcal{C}_{\frac{1}{2^{n}}}(f)\|_{\infty} =\mathcal{O}\left(\frac{1}{2^{2n}}\right).
}

In what follows we will  consider an approach to convolution approximation on sparse grid subsets of  $\mathcal{X}_{n,d}.$ To be more precise, we consider the following subset of $\mathcal{X}_{n,d},$
\EQ{sparse}{
\mathcal{S}_{n,d}=\bigcup_{|\boldsymbol\ell |_{1}=n+d-1}\mathcal{X}_{\boldsymbol\ell},
}
with $|\boldsymbol\ell|_{1}=\ell_{1}+\cdots +\ell_{d},$ which will be referred to as the sparse grid at level $n$ in $d$ dimensions. We note that there is some redundancy in this definition; the sparse grid is represented as a combination of sub-grids and some grid points are included in more than one sub-grid, this is nicely illustrated, for the 2 dimensional case, in Figure 1.

\begin{figure}[h]
\includegraphics[width = 18cm]{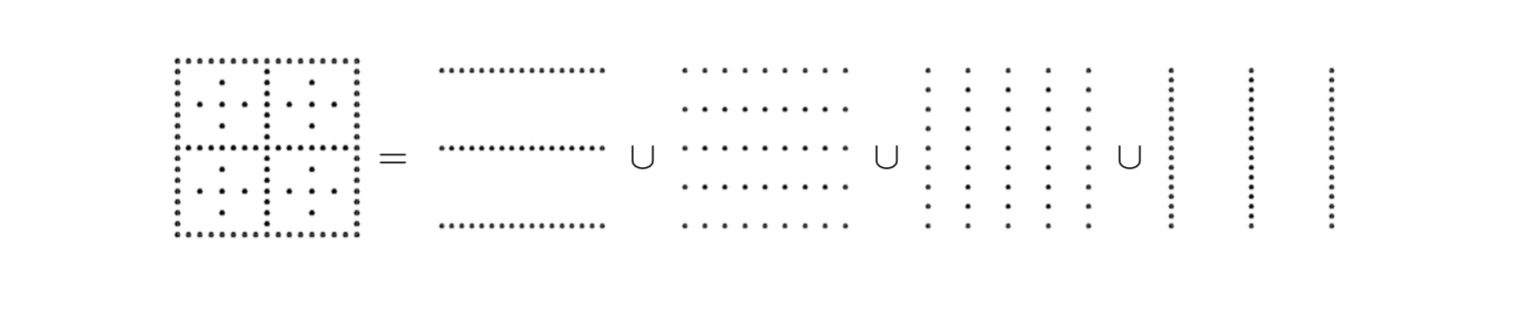}
\centering
\caption{The sparse grid $\mathcal{S}_{4,2}$ constructed via (\ref{sparse})}
\end{figure}

An effort to reduce this redundancy is possible by employing the boolean sum representation of Delvos \cite{delvos}, specifically one can express the sparse grid as

\EQ{sparse-comb}
{
\mathcal{S}_{n,d}=\sum_{q=0}^{d-1}(-1)^{q}{ d-1 \choose q}
\sum_{|\boldsymbol\ell |_{1}=n+(d-1)-q}\mathcal{X}_{\boldsymbol\ell},
}
where we interpret the positive contributions as the inclusion of points and the negative contributions as their removal, this approach is nicely illustrated for the 2 dimensional case, in Figure 2.

\begin{figure}[h]
\includegraphics[width = 18cm]{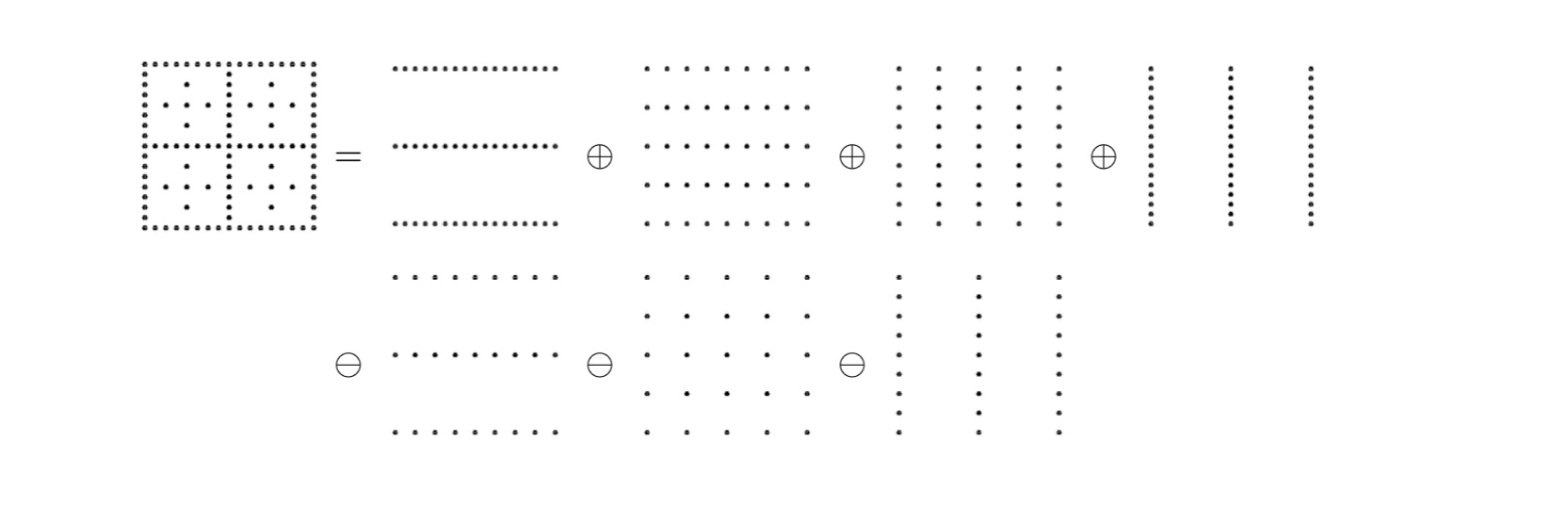}
\centering
\caption{The sparse grid $\mathcal{S}_{4,2}$ constructed via (\ref{sparse-comb})}
\end{figure}

Following (\ref{genconapp}) we represent the anisotropic convolution approximation on $\mathcal{X}_{{\bf{\ell}}}$ as
\EQ{anisotrop}{
\begin{aligned}
\mathcal{C}_{\boldsymbol\ell}(f)(\textbf{x})& = f \ast \Psi_{\left(\frac{1}{2^{\ell_{1}}},\ldots,\frac{1}{2^{\ell_{d}}}\right)}(\textbf{x})=\sum_{\textbf{k}\in {\ZZ}^{d}}\widehat{f}(\textbf{k})\left(\prod_{i=1}^{d}\widehat{\psi}\left(\frac{k_{i}}{2^{\ell_{i}}}\right)
\right)e_{\textbf{k}}(\textbf{x}).
\end{aligned}
}
The convolution approximation on the sparse grid $\mathcal{S}_{n,d}$ makes use of the Boolean decomposition (\ref{sparse-comb}) and, by what is commonly called the {\textit{combination technique}} \cite{griebel}, we define
\EQ{convtec}{
\mathcal{C}_{n,d}(f)(\textbf{x})=\sum_{q=0}^{d-1}(-1)^{q}{ d-1 \choose q}
\sum_{|\boldsymbol\ell |_{1}=n+(d-1)-q}\mathcal{C}_{\boldsymbol\ell}(f)(\textbf{x}).
}
Substituting (\ref{anisotrop}) into the above one can show that
\[
\mathcal{C}_{n,d}(f)(\textbf{x})=\sum_{\textbf{k}\in {\ZZ}^{d}}\widehat{f}(\textbf{k})\widehat{\mathcal{C}_{n,d}}(\textbf{k})e_{\textbf{k}}(\textbf{x})
\]
where
\EQ{combcoef}{
\widehat{\mathcal{C}_{n,d}}(\textbf{k})=(-1)^{d-1}\sum_{q=0}^{d-1}(-1)^{q}{ d-1 \choose q}\sum_{|\boldsymbol\ell|_{1}=n+q}
e^{-\sum_{i=1}^{d}\frac{2\pi^{2}k_{i}^{2}}{2^{2\ell_{i}}}}.
}

Using this representation the pointwise error formula is given by
\EQ{error}{
E_{n,d}(f)({\bf{x}})=f({\bf{x}})-\mathcal{C}_{n,d}(f)(\textbf{x})=\sum_{\textbf{k}\in {\ZZ}^{d}}\widehat{f}(\textbf{k})\widehat{\mathcal{E}_{n,d}}(\textbf{k})e_{\textbf{k}}(\textbf{x}),}
where
\EQ{errorcoef}{
\widehat{\mathcal{E}_{n,d}}(\textbf{k})=1-(-1)^{d-1}\sum_{q=0}^{d-1}(-1)^{q}{ d-1 \choose q}\sum_{|\boldsymbol\ell|_{1}=n+q}
e^{-\sum_{i=1}^{d}\frac{2\pi^{2}k_{i}^{2}}{2^{2\ell_{i}}}}.}
We notice that when $\bf{k}=\bf{0}$ we have that
\EQ{kzero}{\begin{aligned}
\widehat{\mathcal{E}_{n,d}}(\textbf{0})&=1-(-1)^{d-1}\sum_{q=0}^{d-1}(-1)^{q}{ d-1 \choose q}\sum_{|\boldsymbol\ell|_{1}=n+q}1\\&=1-(-1)^{d-1}\sum_{q=0}^{d-1}(-1)^{q}{ d-1 \choose q}{ n+q-1 \choose d-1},
\end{aligned}
}
where we have used the fact that that the number of ways to write $s$ as the sum of $r$ positive integers is ${ s-1 \choose r-1}.$
The following identity, which is taken from \cite{PBMv1}  Formula 4.2.5.47, is valid for non-negative integers $r$ and $s$ such that $0 \le r\le s$
\[
\sum_{q=0}^{s}(-1)^{q}{ s \choose q}{ a+bq \choose r}= (-1)^{s}b^{s}\delta_{r,s}.
\]

Applying this with $r=s = d-1,$ $a = n-1$ and $b=1$ we can conclude that the sum in the expression above equates to $(-1)^{d-1} $ and hence we have that $\widehat{\mathcal{E}_{n,d}}(\textbf{0})=0.$
Thus, as with the plain convolution approximation, the combination convolution approximation on the sparse grid also reproduces the constant function. At this point in the paper it is pertinent to compare the two error representations for convolution approximation that we have developed so far, in the continuous (full grid) setting we have
\[
E_{n,d}(f)({\bf{x}})=f({\bf{x}})-\mathcal{C}_{\frac{1}{2^{n}}}(f)(\textbf{x}) = \sum_{{\bf{k}}\in {\ZZ}^{d} \setminus \{ {\bf{0}}\}}\widehat{f}(\textbf{k})
\left(1-e^{-\frac{2\pi^{2}\|{\bf{k}}\|^{2}}{2^{2n}}}\right)e_{\textbf{k}}(\textbf{x})
\]
and in the sparse grid case we have

\[\begin{aligned}
&E_{n,d}(f)({\bf{x}})=f({\bf{x}})-\mathcal{C}_{n,d}(f)(\textbf{x})=\sum_{{\bf{k}}\in {\ZZ}^{d} \setminus \{ {\bf{0}}\}}\widehat{f}(\textbf{k})\left(1-\sum_{q=0}^{d-1}(-1)^{q+d-1}{ d-1 \choose q}\sum_{|\boldsymbol\ell|_{1}=n+q}
e^{-\sum_{i=1}^{d}\frac{2\pi^{2}k_{i}^{2}}{2^{2\ell_{i}}}}\right)e_{\textbf{k}}(\textbf{x}).
\end{aligned}
\]

In Section 2 we found that error bounds, for sufficiently smooth functions, in the full grid case are easy to access by applying the simple inequality $1-e^{-x}\le x.$ The situation for the sparse grid case is clearly not as straightforward and this leads us to embark on a thorough investigation of the coefficients  (\ref{errorcoef}). To this end we will begin with a detailed examination of the $2-$dimensional case. The findings from the $2-d$ investigation will form the base case of an inductive proof which we will use to establish convergence in higher dimensions.
 

\subsection{Convergence in two-dimensions}

In two dimensions the sparse grid convolution coefficients (\ref{errorcoef}) have the form

\EQ{E2}{
\begin{aligned}
\widehat{\mathcal{E}_{n,2}}(\textbf{k})&=1+\sum_{i+j=n}e^{-2\pi^{2}\left(\frac{k_{1}^{2}}{2^{2i}}+\frac{k_{2}^{2}}{2^{2j}}\right)}-\sum_{i+j=n+1}e^{-2\pi^{2}\left(\frac{k_{1}^{2}}{2^{2i}}+\frac{k_{2}^{2}}{2^{2j}}\right)}.\\
\end{aligned}
}

Let us develop the general term in the above expression using the  full series expansion of the exponential function. Specifically, we consider
\EQ{m2}{
\begin{aligned}
\sum_{i+j=m}e^{-2\pi^{2}\left(\frac{k_{1}^{2}}{2^{2i}}+\frac{k_{2}^{2}}{2^{2j}}\right)}&=
\sum_{i+j=m}\sum_{p=0}^{\infty}(-1)^{p}\frac{2^{p}\pi^{2p}}{p!}\left(\frac{k_{1}^{2}}{2^{2i}}+\frac{k_{2}^{2}}{2^{2j}}\right)^{p}=
\sum_{p=0}^{\infty}(-1)^{p}\frac{2^{p}\pi^{2p}}{p!}\sigma_{2,p}(m,\textbf{k}),
\end{aligned}
}
where
\EQ{sig2}{
\sigma_{2,p}(m,\textbf{k})=\sum_{i+j=m}\left(\frac{k_{1}^{2}}{2^{2i}}+\frac{k_{2}^{2}}{2^{2j}}\right)^{p}=\sum_{j=1}^{m-1}\left(\frac{k_{1}^{2}}{2^{2(m-j)}}+\frac{k_{2}^{2}}{2^{2j}}\right)^{p}.
}
We note that, when $p=0$ we have $\sigma_{2,0}(m,\textbf{k})=m-1.$ For $p\ge 1$ we can apply the binomial theorem to yield
\[
\sigma_{2,p}(m,\textbf{k})=\sum_{j=1}^{m-1}\sum_{r=0}^{p}{ p \choose r}\frac{k_{1}^{2r}}{2^{2(m-j)r}}\frac{k_{2}^{2(p-r)}}{2^{2j(p-r)}}
=\sum_{r=0}^{p}{ p \choose r}\frac{k_{1}^{2r}k_{2}^{2(p-r)}}{2^{2mr}}\sum_{j=1}^{m-1}2^{2(2r-p)j}.
\]
Define
\EQ{delta}{\Delta_{p,2}=\begin{cases} 1\,\, & \,\,
\textrm{if} \,\,p\,\,{\rm{is}}\,\,{\rm{even}}; \\
0 \,\, & \,\,
\textrm{otherwise},
\end{cases}}
then applying the geometric sum formula
\EQ{gsum}{
\sum_{j=1}^{n}x^{j}=\begin{cases} \frac{1}{x^{-1}-1}(1-x^{n})\,\, & \,\,
\textrm{if} \,\,x\ne 1; \\
n \,\, & \,\,\textrm{if} \,\,x= 1.
\end{cases}}
We find that
\[
\begin{aligned}
&\sigma_{2,p}(m,\textbf{k})=
\sum_{\substack{r=0\\ r\ne \frac{p}{2}}}^{p}{ p \choose r}\frac{k_{1}^{2r}k_{2}^{2(p-r)}}{2^{2mr}}\frac{1-2^{2(2r-p)(m-1)}}{2^{2(p-2r)}-1}+\Delta_{p,2}{ p \choose \frac{p}{2}}\frac{k_{1}^{p}k_{2}^{p}}{2^{mp}}(m-1)\\
&=\sum_{\substack{r=0\\ r\ne \frac{p}{2}}}^{p}{ p \choose r}\left(\frac{k_{1}^{2r}k_{2}^{2(p-r)}}{2^{2mr}(2^{2(p-2r)}-1)}-
\frac{k_{1}^{2(p-r)}k_{2}^{2r}2^{2(p-2r)(m-1)}}{2^{2m(p-r)}(2^{2(2r-p)}-1)}\right)
+\Delta_{p,2}{ p \choose \frac{p}{2}}\frac{k_{1}^{p}k_{2}^{p}}{2^{mp}}(m-1)\\
&=\sum_{\substack{r=0\\ r\ne \frac{p}{2}}}^{p}{ p \choose r}\frac{k_{1}^{2r}k_{2}^{2(p-r)}+k_{1}^{2(p-r)}k_{2}^{2r}}{2^{2mr}(2^{2(p-2r)}-1)}+
\Delta_{p,2}{ p \choose \frac{p}{2}}\frac{k_{1}^{p}k_{2}^{p}}{2^{mp}}(m-1).
\end{aligned}
\]

Using the notation introduced above we can write (\ref{E2}) as
\[
\begin{aligned}
\widehat{\mathcal{E}_{n,2}}(\textbf{k})=&\sum_{p=1}^{\infty}\frac{(-1)^{p}2^{p}\pi^{2p}}{p!}\left(\sigma_{2,p}(n,\textbf{k})-\sigma_{2,p}(n+1,\textbf{k})\right)\\
=&\sum_{p=1}^{\infty}\frac{(-1)^{p}2^{p}\pi^{2p}}{p!}\sum_{\substack{r=1\\ r\ne \frac{p}{2}}}^{p}\frac{{ p \choose r}}{2^{2nr}}\left(\frac{2^{2r}-1}{2^{2(p-r)}-2^{r}}\right)\left(k_{1}^{2r}k_{2}^{2(p-r)}+k_{1}^{2(p-r)}k_{2}^{2r}\right)\\
&+\sum_{p=1}^{\infty}\frac{2^{2p}\pi^{4p}}{(2p)!}{ 2p \choose p}\left(\left(\frac{2^{2p}-1}{2^{2p}}\right)\frac{n}{2^{2np}}-\frac{1}{2^{2np}}\right)k_{1}^{2p}k_{2}^{2p}.
\end{aligned}
\]

Examining the above, we observe that the term  dominating the asymptotic rate of decay corresponds to the first ($p=1$) term of the second sum, and hence we can deduce that
\[
\widehat{\mathcal{E}_{n,2}}(\textbf{k})\sim \frac{ 2^{2}\pi^{4}}{2!}2\cdot \frac{3}{4}\frac{n}{2^{2n}}k_{1}^{2}k_{2}^{2}=3\pi^{4}k_{1}^{2}k_{2}^{2}\frac{n}{2^{2n}}.
\]


\subsection{Convergence in $d$-dimensions}

In this part we will consider the $d-$dimensional analogue of the approach from the previous subsection. We begin by defining the $d-$dimensional analogue of (\ref{sig2})
\EQ{sigd}{
\sigma_{d,p}(m,\textbf{k})=\begin{cases}{ m-1 \choose d-1}\,\, & \,\,
\textrm{if} \,\,p=0; \\
\sum_{|\boldsymbol\ell|_{1}=m}\left(\sum_{i=1}^{d}\frac{k_{i}^{2}}{2^{2\ell_{i}}}\right)^{p}\,\, & \,\,
\textrm{if} \,\,p\ge 1.
\end{cases}}

Then, using this notation in the expansion of the exponential function, the error coefficients (\ref{errorcoef}) can be represented as
\[
\begin{aligned}
\mathcal{E}_{n,d}(\textbf{k})=&1-(-1)^{d-1}\sum_{p=0}^{\infty}\frac{(-1)^{p}2^{p}\pi^{2p}}{p!}\sum_{q=0}^{d-1}(-1)^{q}{ d-1 \choose q}\sigma_{d,p}(n+q,\textbf{k})\\
=&1-(-1)^{d-1}\sum_{q=0}^{d-1}(-1)^{q}{ d-1 \choose q}{ n+q-1 \choose d-1}\\
&-\sum_{p=1}^{\infty}\frac{(-1)^{p}2^{p}\pi^{2p}}{p!}\sum_{q=0}^{d-1}(-1)^{d-1-q}{ d-1 \choose q}\sigma_{d,p}(n+q,\textbf{k}).
\end{aligned}
\]
We note that the penultimate line above coincides with $\mathcal{E}_{n,d}(\textbf{0})$ (\ref{kzero}) which we have shown to be zero.
To simplify the notation we recall the  forward divided difference functional of order $k$ is defined by
\[\Delta^{k}f= \sum_{q=0}^{k}(-1)^{k-q}{ k \choose q}f(q).
\]
Taking $k=d-1$ we can express the $d-$dimensional sparse grid convolution error coefficients as
\[
\mathcal{E}_{n,d}(\textbf{k})=-2\pi^{2}\Delta^{d-1}\sigma_{d,1}(n+\cdot,\textbf{k})-\sum_{p=2}^{\infty}\frac{(-1)^{p}2^{p}\pi^{2p}}{p!}\Delta^{d-1}\sigma_{d,p}(n+\cdot,\textbf{k}).
\]
Clearly, an investigation of $\sigma_{d,p}(m,\textbf{k})$ is required in order to shed light upon the rate at which the error coefficients decay. The following result provides the insight we need.

\

\theo{Theorem}{main1}{
Let $d\ge 2,$  $\boldsymbol\ell=(\ell_{1},\ldots,\ell_{d})^{T} \in \N^{d}$ and
${\textbf{k}}=(k_{1},\ldots,k_{d})^{T} \in {\ZZ}^{d}.$ Then, for $m\ge d,$
we have
\EQ{linbit}{
\sigma_{d,1}(m,\textbf{k})
=\|\textbf{k}\|^{2}\left(\pi_{d-2}(m)+(-1)^{d-1}\left(\frac{4}{3}\right)^{d-1}\frac{1}{2^{2m}}\right),}
where $\pi_{d-2}(m)$ is a polynomial in $m$ of degree $d-2.$  Furthermore, for $p\ge 2$ we also have
\EQ{rest}{
\sigma_{d,p}(m,\textbf{k})=\pi_{d-2}^{\textbf{k},p}(m)+\delta_{p,d} \frac{d k_{1}^{2}\cdots k_{d}^{2}\, m^{d-1}}{2^{2m}}+
 \frac{C(\textbf{k},p,d)m^{d-2}}{2^{2m}}+O\left(\frac{m^{d-3}}{2^{2m}}\right),
}
where $\pi_{d-2}^{\textbf{k},p}(m)$ is a polynomial in $m$ of degree $d-2$ whose coefficients depend upon $\textbf{k}$ and $p.$}

The full proof of this theorem relies on some rather technical machinery and this, together with the proof, is provided in the final section of the paper. Using the expression for $\sigma_{2,p}(m,\textbf{k})$ that was derived in the previous subsection it is easy to verify the $2-$dimensional version of the result with constant polynomials
\[
\pi_{0}(m) =\frac13,\,\,{\rm{and}}\,\,\pi_{0}^{\textbf{k},p}(m)= \frac{k_{1}^{2p}+k_{2}^{2p}}{2^{2p}-1},\,\,\, p\ge 2,
\]
and leading constants
\[C(\textbf{k},p,2)=
\frac{p\left(k_{1}^{2(p-1)}k_{2}^{2}+k_{1}^{2}k_{2}^{2(p-1)}\right)}{2^{2(p-2)}-1},\,\, p\ne 2.
\]
The key insight from Theorem \ref{main1} is  that, for each $p,$ the $\sigma_{d,p}(m,\textbf{k})$ function can be expressed as a polynomial in $m$ of degree $d-2$ plus either a constant multiple of $m^{d-2}/2^{2m}$ (when $p\ne d$) or $dk_{1}^{2}\cdots k_{d}^{2} m^{d-1}/2^{2m}$ (when $p=d)$ followed by higher order terms (i.e., those decaying at a faster rate as $m$ grows).  Given that the forward divided difference functional annihilates polynomials of degree $d-2$ we can, after ignoring the higher order terms, deduce that
\[
\begin{aligned}
\widehat{\mathcal{E}_{n,d}}(\textbf{k})=\frac{1}{2^{2n}}&\Biggl[(-1)^{d}2\pi^{2}\|\textbf{k}\|^{2}\sum_{q=0}^{d-1}(-1)^{d-1-q}{ d-1 \choose q}\frac{1}{2^{2q}}\\
&-\sum_{\substack{p=2\\ p\ne d}}^{\infty}\frac{(-1)^{p}2^{p}\pi^{2p}C(\textbf{k},p,d)}{p!}\sum_{q=0}^{d-1}(-1)^{d-1-q}{ d-1 \choose q}\frac{(n+q)^{d-2}}{2^{2q}}\\&-\frac{(-1)^{d}2^{d}\pi^{2d}dk_{1}^{2}\cdots k_{d}^{2}}{d!}\sum_{q=0}^{d-1}(-1)^{d-1-q}{ d-1 \choose q}\frac{(n+q)^{d-1}}{2^{2q}}+O(n^{d-3})\Biggr].
\end{aligned}
\]

Examining the above sum we see that the asymptotic decay of the coefficients is dominated by $n^{d-1}$ weight arising from  the application of the forward divided difference operator  to $(n+\cdot\,)^{d-1}$ in the final sum. Thus, using the binomial identity
\[
\sum_{q=0}^{d-1}(-1)^{d-1-q}{ d-1 \choose q}x^{q}= (-1)^{d-1}(1-x)^{d-1},
\] with $x=1/4$ we may deduce that
\EQ{ass}{
\widehat{\mathcal{E}_{n,d}}(\textbf{k})\sim \frac{2\pi^{2d}k_{1}^{2}\cdots k_{d}^{2}}{(d-1)!}\left(\frac{3}{2}\right)^{d-1}\frac{n^{d-1}}{2^{2n}}.
}

Employing this result in (\ref{error})  we can deduce

\EQ{errass}{
E_{n,d}(f)({\bf{x}})\sim \frac{2\pi^{2d}}{(d-1)!}\left(\frac{3}{2}\right)^{d-1}\frac{n^{d-1}}{2^{2n}}\sum_{\textbf{k}\in Z^{d}}\widehat{f}(\textbf{k})k_{1}^{2}\cdots k_{d}^{2}e_{\textbf{k}}(\textbf{x}),}
and,  more specifically, by mirroring the proof of Proposition \ref{prop1}, we can conclude the following.

\

\theo{Corollary}{cor1}{
Let $f \in \mathcal{N}_\beta,$ where $\beta >\frac{d}{2}+2.$ Let
$\mathcal{C}_{\frac{1}{2^{n}}}(f)$ denote the plain Gaussian convolution approximation (\ref{genconapp}) to $f$ on the full isotropic grid   $\mathcal{X}_{n,d}$ with spacing  $1/2^{n}$  and  $\mathcal{C}_{n,d}(f)$  denote the combined convolution approximation to $f$ (\ref{convtec}) on the sparse grid $\mathcal{S}_{n,d}.$ Then
\[
\|\mathcal{C}_{\frac{1}{2^{n}}}(f)-f\|_{\infty}\le \frac{C_{d}}{2^{2n}}\|f\|_{\beta}\quad {\rm{and}}\quad\|\mathcal{C}_{n,d}(f)-f\|_{\infty} \le \frac{C_{d}n^{d-1}}{2^{2n}}\|f\|_{\beta},
\]
where $C_{d}$ denotes a generic dimension dependent constant.
}

\

We close this section by presenting some numerical results to show  how closely the Fourier coefficients of the sparse grid convolution approximation track the asymptotic formula.

\begin{table}[h]
\caption{Comparison of numerically computed $2-$d expansion coefficients  $\widehat{\mathcal{E}_{n,2}}(\textbf{k})$ with the  asymptotic formula  (\ref{ass}), with $\textbf{k}=(1,1)$ (left) and $\textbf{k}=(500,700)$ (right)}
\centering
\begin{tabular}{|c|c|c|c|c|} \hline
 $n$  & $\widehat{\mathcal{E}_{n,2}}(\textbf{k})$ &  formula & $\widehat{\mathcal{E}_{n,2}}(\textbf{k})$ & formula  \\ \hline \hline
40 & 8.69 (-21) & 9.67 (-21)  & 5.19 (-10) & 1.18 (-9)  \\ \hline
80 & 1.52 (-44) & 1.60 (-44) & 1.41 (-33) & 1.99 (-33) \\ \hline
160 & 2.13 (-92) & 2.19 (-92)& 2.31 (-81) & 2.68 (-81)  \\ \hline
320 & 2.02 (-188) & 2.05 (-188) & 2.33 (-177) & 2.51 (-177) \\ \hline
640 & 8.93 (-381) & 8.98 (-381) & 1.06 (-369) & 1.10 (-369) \\ \hline
\end{tabular}
\end{table}

\begin{table}[h]
\caption{Comparison of numerically computed $3-$d expansion coefficients  $\widehat{\mathcal{E}_{n,3}}(\textbf{k})$ with the  asymptotic formula  (\ref{ass}), with $\textbf{k}=(1,1,1)$ (left) and $\textbf{k}=(500,700,900)$ (right)}
\centering
\begin{tabular}{|c|c|c|c|c|} \hline
 $n$  & $\widehat{\mathcal{E}_{n,3}}(\textbf{k})$ &  formula & $\widehat{\mathcal{E}_{n,3}}(\textbf{k})$ & formula  \\ \hline \hline
40 & 1.72 (-18) & 2.86 (-18)  & 3.54 (-2) & 2.84 (-1)  \\ \hline
80 & 6.65 (-42) & 9.47 (-42) & 4.62 (-25) & 9.40 (-25) \\ \hline
160 & 1.95 (-89) & 2.59 (-89)&1.69 (-72) & 2.57 (-72)  \\ \hline
320 & 3.80 (-185) & 4.85 (-185) &3.54 (-168) & 4.82 (-168) \\ \hline
640 & 3.39 (-377) & 4.26 (-377) &3.27 (-360) & 4.22 (-360) \\ \hline
\end{tabular}
\end{table}

\newsection{Proof of Main Theorem}

The main results stated in Theorem \ref{main1} are not so hard to convey.  The polynomials that appear in the results arise from terms in the multinomial expansion which involve iterations of finite geometric series; some of these collapse to the sum of powers of natural numbers and, as such, introduce polynomial terms in $m.$  For instance, in the $2-$dimensional investigated in subsection 3.2, we see the finite geometric sum $\sum_{j=1}^{m-1}2^{2(2r-p)}$ and, in the cases where $2r=p$ this collapses to $\sum_{j=1}^{m-1}1=m-1,$ thus introducing a linear term in $m.$ If one was to carefully examine the $3-d$ case then further instances of such sums (leading to linear terms in $m$) would arise together with double sums that collapse to  $\sum_{j=1}^{m-2}j=(m-1)(m-2)/2,$ and these introduce a quadratic term. The pattern continues into higher dimensions.  The formal proof of the result is made difficult due, in part, to the notational complexity that is involved. The first result, identity (\ref{linbit}), is a surprisingly neat representation for the $p=1$ case; we were not able to develop similar neat closed form expressions for $p>2.$ We begin by establishing (\ref{linbit}), then we will develop some technical results on weighted geometric sums that will allow us to verify (\ref{rest}). We begin with the following lemma which sheds some insight on a particular finite sum.

\

\theo{Lemma}{l1}
{Let $d\ge2$ be a positive integer and $m> d. $ Then, for a positive integer $r,$ we have
\[
\sum_{j=1}^{m-(d-1)}\frac{{ m-j-1 \choose d-2}}{2^{2jr}}=p_{d-2}(m)+\left(\frac{2^{2r}}{1-2^{2r}}\right)^{d-1}\frac{1}{2^{2rm}},\]
where $p_{d-2}(m)$ is a polynomial in $m$ of degree $d-2.$}
\begin{proof}
Recall  the Gauss
hypergeometric
    function (see \cite{AS}, 15.1.1) is
    defined by
    \EQ{hyper}{
_{2}F_{1}(a,b;c;z):=\sum_{j=0}^{\infty}\frac{
(a)_j (b)_{j} }{(c)_j}\frac{z^{j}}{j!}, }  where \EQ{poch}{ (x)_j := x(x+1)\cdots(x+j-1) \quad
j \geq 1} denotes the
Pochhammer symbol, with $(x)_0 = 1$. If $n$ is a positive integer we have
\[
(n)_{j}=\frac{(n+j-1)!}{(n-1)!}\quad{\rm{and}}\quad (-n)_{j}=\frac{(-1)^{j}n!}{(n-j)!}.
\]
Using the above it is straight forward to verify that
\EQ{hyp}{
\sum_{j=1}^{m-(d-1)}\frac{{ m-j-1 \choose d-2}}{2^{2jp}}=\frac{1}{2^{2p}}{ m-2 \choose d-2}{}_2F_1\left(1,-(m-d);-(m-2); \frac{1}{2^{2p}}\right).
}

The following identity, see \cite{PBMv3} Formula 7.3.1.178, is valid for non-negative integers $p$ and $q$
\[
_{2}F_{1}(1,-p;-q;z)=\frac{q+1}{p+1}\sum_{k=0}^{q-p}\frac{(p-q)_{k}}{(p+2)_{k}}(1-z)^{-(k+1)}+\frac{(-1)^{p}p!}{(-q)_{p}}z^{q+1}(z-1)^{p-q-1}.
\]

Applying this with $p=m-d,$ $q=m-2$ and $z=\frac{1}{2^{2r}},$ we find that
\[
\begin{aligned}
&_{2}F_1\left(1,-(m-d);-(m-2); \frac{1}{2^{2r}}\right)=\\&\frac{m-1}{m-(d-1)}\sum_{k=0}^{d-2}\frac{\left(-(d-2)\right)_{k}\left(\frac{2^{2r}}{2^{2r}-1}\right)^{k+1}}{\left(m-(d-2)\right)_{k}}+
\frac{(-1)^{m-1}(m-d)!\left(\frac{2^{2r}}{2^{2r}-1}\right)^{d-1}}{\left(-(m-2)\right)_{m-d}}\frac{2^{2r}}{2^{2rm}}\\
&=(m-1)\sum_{k=0}^{d-2}\frac{(-1)^{k}(d-2)!(m-d)!\left(\frac{2^{2r}}{2^{2r}-1}\right)^{k+1}}{(d-2-k)!(m+k-(d-1))!}
+\frac{(-1)^{d-1}2^{2r}\left(\frac{2^{2r}}{2^{2r}-1}\right)^{d-1}}{{ m-2 \choose d-2}2^{2rm}}.
\end{aligned}
\]

In view of (\ref{hyp}) we now multiply this by $\frac{1}{2^{2r}}{ m-2 \choose d-2}$ and, following some elementary
 simplifications, we have  the following expression
 \[
 \begin{aligned}
 \sum_{j=1}^{m-(d-1)}\frac{{ m-j-1 \choose d-2}}{2^{2rj}}&=
 \sum_{k=0}^{d-2}\frac{(-1)^{k}{ m-1 \choose m+k-(d-1)}\left(\frac{2^{2r}}{2^{2r}-1}\right)^{k}}{2^{2r}-1}
+\frac{(-1)^{d-1}\left(\frac{2^{2r}}{2^{2r}-1}\right)^{d-1}}{2^{2rm}}\\
 & = p^r_{d-2}(m)+\left(\frac{2^{2r}}{1-2^{2r}}\right)^{d-1}\frac{1}{2^{2rm}},
 \end{aligned}
 \]
where  $p^r_{d-2}(m),$ which represents the sum appearing above, is a polynomial in $m$ of degree $d-2.$
 \end{proof}

\subsection{Proof of (\ref{linbit})}
We know from $(3.30)$ that
\[
\begin{aligned}
&\sigma_{d,1}(m,\textbf{k})=\sum_{|\boldsymbol\ell|_{1}=m}\frac{k_{1}^{2}}{2^{2\ell_{1}}}+\cdots +\frac{k_{d}^{2}}{2^{2\ell_{d}}}=\sum_{i=1}^{d}k_i^2\sum_{|\boldsymbol\ell|_{1}=m}\frac{1}{2^{2\ell_i}.}
\end{aligned}\]

The above sum concerns the set of $d-$dimensional multi-indices $\boldsymbol\ell$ satisfying $|\boldsymbol\ell|_{1}=m.$ A typical component $\ell_{i}$ of  $\boldsymbol\ell$ can, theoretically, take on any value between $1$ and $m-(d-1)$ included (in the latter case remaining  $d-1$ components are all set to $1)$. The number of times $\ell_{i}$ takes on a certain value $j \in \{1,2,\ldots,m-(d-1)\}$ is precisely the number of ways in which the remaining $d-1$ components of  $\boldsymbol\ell$ sum to $m-j$ and this is given by $
{m-j-1 \choose d-2}.$ Since the last sum in the above expression only depends on the value $\ell_i$ and not on $i$ we have that
\[
\begin{aligned}
\sigma_{d,1}(m,\textbf{k})
&=\sum_{i=1}^{d}k_i^2\sum_{j=1}^{m-(d-1)}{m-j-1 \choose d-2} \frac{1}{2^{2j}}=\|\textbf{k}\|^{2}\left( p_{d-2}(m)+\left( \frac{2^2}{1-2^2}\right)^{d-1}\frac{1}{2^{2m}}\right),
\end{aligned}
\]
where the last equation follows from Lemma 4.1, with $r=1,$ and the proof of (\ref{linbit}) is complete.
 \subsection{On weighted geometric sums}

In this subsection we outline some key results on the representation of the kinds of weighted geometric series that are encountered if one applies the appropriate multinomial expansion in order to examine the sums $\sigma_{d,p}(m,\bf{k})$ (\ref{sigd}) for $p\ge 2.$ We begin by differentiating the plain geometric sum formula, followed by multiplication by $x$ to deduce that

 \EQ{firstcase}{
\sum_{j=1}^{n}jx^{j}=x\frac{d}{dx}\sum_{j=1}^{n}x^{j}=\begin{cases} \frac{x}{(1-x)^{2}}\left(1-x^{n}\left((1+n)-nx\right)\right) \,\, & \,\,
{\textrm{if}} \,\,\,\, x\ne 1; \\
\sum_{k=1}^{n}k = \frac{n(n+1)}{2} \,\, & \,\,
{\textrm{if}} \,\,\,\, x=1.
\end{cases} 
}

Let us consider the more general weighted geometric sum
\[
G_{i}^{(n)}(x)=\sum_{j=1}^{n}j^{i}x^{j}.
\]

We note in the case where $x=1$ we have the sum of the $i^{th}$ powers of the first $n$ positive integers which, due to Faulhaber's formula, see \cite{GR65} Formula 0.121, is a polynomial in $n$ of degree $i+1,$

\EQ{Faul}{
G_{i}^{(n)}(1)=\sum_{j=1}^{n}j^{i}=\frac{n^{i+1}}{i+1}+q_{i}(n),
}
where $q_{i}(n)$ is a polynomial in $n$ of degree $i.$
For the more general case ($x \ne 1$) we observe that
\EQ{drel}{
G_{i+1}^{(n)}(x)=x\frac{d}{dx}G_{i}^{(n)}(x)
} and this allows us to establish the following.
\theo{Lemma}{sums}{Let $j$ denote a non-negative integer and $x\ne 1,$ then
\EQ{rep}{G_{i}^{(n)}(x)=\frac{x}{(1-x)^{i+1}}\left(q_{i-1}(x)-x^{n}p_{i}(x,n)\right),}
where $q_{i-1}(x)$ is a polynomial of degree $i-1$  in $x$  and $p_{i}(x,n)$ is a polynomial of degree $i$ in both $x$ and $n.$
}

\begin{proof}
We establish the result via induction. Appealing to (\ref{firstcase}) we see that the result is true for $i=1$ with $q_{0}(x)=1$ and
$p_{1}(x,n)=1+n-nx.$ Assume the result is true for $i$ and consider the following development, using (\ref{drel}), for the case $i+1.$
\[
\begin{aligned}
G_{i+1}^{(n)}(x)&=x\frac{d}{dx}G_{i}^{(n)}(x)=x\frac{d}{dx}\left(\frac{x}{(1-x)^{i+1}}\left(q_{i-1}(x)-x^{n}p_{i}(x,n)\right)\right)\\
&=x\Bigl[\frac{x}{(1-x)^{i+1}}\left(\frac{d}{dx}q_{i-1}(x)-x^{n}\frac{d}{dx}p_{i}(x,n)-nx^{n-1}p_{i}(x,n)\right)\\
&+\left(\frac{1}{(1-x)^{i+1}}+\frac{ix+x}{(1-x)^{i+2}}\right)\left(q_{i-1}(x)-x^{n}p_{i}(x,n)\right)\Bigr]\\
&=\frac{x}{(1-x)^{i+2}}\Bigl[x(1-x)\left(\frac{d}{dx}q_{i-1}(x)-x^{n}\frac{d}{dx}p_{i}(x,n)-nx^{n-1}p_{i}(x,n)\right)\\
&+(1+ix)\left(q_{i-1}(x)-x^{n}p_{i}(x,n)\right)\Bigr]\\
&=\frac{x}{(1-x)^{i+2}}\Bigl[q_{i}(x)-x^{n}p_{i+1}(x,n)\Bigr],
\end{aligned}
\]
where
\[
q_{i}(x)=(1+ix)q_{i-1}(x)+x(1-x)\frac{d}{dx}q_{i-1}(x)
\] is clearly a polynomial in $x$ of degree $i$ and, likewise, where
\[
p_{i+1}(x,n)=(1+ix +n-nx)p_{i}(x,n)+x(1-x)\frac{d}{dx}p_{i}(x,n)
\]is clearly a polynomial of degree $i+1$ in both $x$ and $n.$

\end{proof}

In order to prepare for how the above result will be used, we let $d$ be the fixed spatial dimension and $m>d$ a positive integer. In what follows we will evaluate various sums and, in each case, we will ignore terms that decay faster than $m^{d-2}2^{-2m}$ for large $m.$  In each case we consider a fixed integer parameter $t\ge 1$ and, where appropriate, we will also consider specific cases of $t=0$ and $t=-1.$ We begin with a straightforward geometric sum for $t\geq1$
\EQ{v1}{
\sum_{j=1}^{m-d}\frac{1}{2^{2jt}}=\frac{1}{2^{2t}-1}\left(1-\frac{2^{2dt}}{2^{2mt}}\right) = {\rm{constant}}+O(2^{-2m}).
}

For the following sum with $t>1$ we can directly use (\ref{rep}) in its evaluation:
\EQ{v2}{
\begin{aligned}
\sum_{j=1}^{m-d}\frac{(m-j)^{i}}{2^{2tj}}&=\frac{1}{2^{2mt}}\sum_{j=d}^{m-1}2^{2tj}j^{i}=\frac{1}{2^{2mt}}\left(G_{i}^{(m-1)}(2^{2t})-G_{i}^{(d-1)}(2^{2t})\right)\\
&=\frac{1}{(1-2^{2t})^{i+1}}\Bigl[\frac{2^{2td}}{2^{2tm}}p_{i}(2^{2t},d-1)-p_{i}(2^{2t},m-1)\Bigr]\\
&=\frac{P_{i}(m)}{(1-2^{2t})^{i+1}}+O(2^{-2m}),
\end{aligned}
}
where $P_{i}(m)$ is a polynomial in $m$ of degree $i$ whose coefficients depend on $2^{2t}.$ In the case where $t=0$ the above collapses to
\EQ{v3}{
\begin{aligned}
\sum_{j=1}^{m-d}(m-j)^{i}&=\sum_{j=d}^{m-d}j^{i}=G_i^{m-1}(1)-G_i^{d-1}(1)\\
&=\frac{(m-1)^{i+1}}{i+1}-\frac{(d-1)^{i+1}}{i+1}+q_{i}(0,d-1)-q_{i}(0,m-1)\\&
=\frac{m^{i+1}}{i+1}+\pi_i(m),
\end{aligned}
}
where $\pi_{i}(m)$ is a polynomial in $m$ of degree $i.$  In the case where $t=-1$ we have
\EQ{v4}
{
\sum_{j=1}^{m-d}2^{2j}(m-j)^{i}=C_{d,i}2^{2m}+\left(\frac{4}{3}\right)^{i+1}P^{*}_{i}(m)}
where $P^{*}_{i}(m)$ is a polynomial in $m$ of degree $i$ and $C_{d,i}$ is a constant depending on $d$ and $i.$\\

\subsection{Proof of (\ref{rest})}
We know from subsection 3.1. that  (\ref{rest}) holds for $d=2$ and any value $p$, let us also assume that it is true for $d_{0}<d$ and any value of $p,$ we will now proceed to show, by induction, that the same statement  is true for $d$ and any value of $p.$ First we establish a recurrence relation for $\sigma_{d,p}$ using $\boldsymbol r=(r_1,\ldots,r_d)^T$,  $\hat{\boldsymbol r}=(r_1,\ldots,r_{d-1})^T$, $\hat{\boldsymbol\ell}=(\ell_1,\ldots,\ell_{d-1})^T$, $\hat{\boldsymbol k}=(k_1,\ldots,k_{d-1})^T$ and  the  multinomial theorem to find
 \[
\begin{aligned}
&\sigma_{d,p}(m,\textbf{k})=\sum_{|\boldsymbol\ell|_{1}=m}\left(\sum_{i=1}^{d}\frac{k_{i}^{2}}{2^{2\ell_{i}}}\right)^{p}\\
&=\sum_{|\boldsymbol\ell|_{1}=m} \sum_{|r|=p}{p \choose r_1\ \ldots\ r_d} \prod_{i=1}^d \left( \frac{k_i^2}{2^{2\ell_i}}\right)^{r_i}\quad \text{( multinomial theorem)}\\
&= \sum_{|\boldsymbol r|=p}{p \choose r_1\ \ldots\ r_d} \prod_{i=1}^d k_i^{2r_i} \sum_{|\boldsymbol\ell|_{1}=m} \left( 2^{\sum_{i=1}^d -2\ell_ir_i}\right)\\
&= \sum_{r_d=0}^{p} \frac{p! k_d^{2r_d}}{(p-r_d)!r_d!}\sum_{\ell_d=1}^{m-d}2^{-2\ell_d r_d} \sum_{|\hat{\boldsymbol r}|=p-r_d} {p- r_d \choose r_1\ \ldots\ r_{d-1}} \prod_{i=1}^{d-1} k_i^{2r_i}  \sum_{|\hat{\boldsymbol \ell}|_{1}=m-\ell_d} \left( 2^{\sum_{i=1}^d -2\ell_ir_i}\right) \\
&= \sum_{r_d=0}^{p}   \frac{p! k_d^{2r_d}}{(p-r_d)!r_d!} \sum_{\ell_d=1}^{m-d}\frac{\sigma_{d-1,p-r_d}(m-\ell_d,\hat{\boldsymbol k})}{2^{2\ell_d r_d}}.
\end{aligned}
\]
Applying the inductive hypothesis (\ref{rest}) we have
\[\begin{aligned}
 \sigma_{d-1,p-r_d}(m-\ell_d,\hat{\boldsymbol k})=&  \pi_{d-3}^{\hat{\boldsymbol k},p-r_d}(m-\ell_d)+\delta_{p-r_d,d-1}\frac{(d-1)k_{1}^2\cdots k_{d-1}^2(m-\ell_d)^{d-2}}{2^{2(m-\ell_d)}}\\
 &+  C(\hat{\boldsymbol k},p-r_d,d-1)\frac{(m-\ell_d)^{d-3}}{2^{2(m-\ell_d)}} +O\left(\frac{m^{d-4}}{2^{2m}}\right).
\end{aligned}\]
Inserting this representation into the inner sum of the above computation yields
\EQ{startpoint}{
\sum_{\ell_d=1}^{m-d}\frac{\sigma_{d-1,p-r_d}(m-\ell_d,\hat{\boldsymbol k})}{2^{2\ell_d r_d} }=S_{1}(r_{d})+S_{2}(r_{d})+S_{3}(r_{d})+O\left(\frac{m^{d-3}}{2^{2m}}\right),}
where
\EQ{S123}{
\begin{aligned}
S_{1}(r_{d}) &=   \sum_{\ell_d=1}^{m-d}  \frac  {\pi_{d-3}^{\hat{\boldsymbol k},p-r_d}(m-\ell_d)}{2^{2\ell_d r_d}},\\
S_{2}(r_{d})&=\delta_{p-r_d,d-1} (d-1)\frac{k_{1}^2\cdots k_{d-1}^2}{2^{2m}}\sum_{\ell_d=1}^{m-d}\frac{(m-\ell_d)^{d-2}}{2^{2\ell_d(r_{d}-1)}},\\
{\rm{and}}\quad S_{3}(r_{d})&=\frac{C(\hat{\boldsymbol k},p-r_d,d-1)}{2^{2m}}\sum_{\ell_d=1}^{m-d}\frac{(m-\ell_d)^{d-3}}{2^{2\ell_d (r_{d}-1)}}. 
\end{aligned}
}

For  the final  term of (\ref{startpoint}) we have used that the sum consists of less than $m$ terms and  each of which is order $\frac{m^{d-4}}{2^{2m}}$. We can use the weighted geometric sums to investigate the three summands above. We begin with $S_{1}(r_{d}) $ and in this case we write the polynomial of degree $d-3$ as
\[
\pi_{d-3}^{\hat{\boldsymbol k},p-r_d}(m-\ell_d)=\sum_{i=0}^{d-3}a_{i}^{\hat{\boldsymbol k},p-r_d}(m-\ell_d)^{i},
\]
and thus we have
\[
S_{1}(r_{d})=\sum_{i=1}^{d-3} a_{i}^{\hat{\boldsymbol k},p-r_d}\sum_{\ell_d=1}^{m-d} \frac{(m-\ell_d)^{i}}{ 2^{2\ell_dr_d}}.
\]
Appealing to (\ref{v2}) and (\ref{v3}) we have that
\[
\sum_{\ell_d=1}^{m-d} \frac{(m-\ell_d)^{i}}{ 2^{2\ell_dr_d}}=\begin{cases}
\frac{P_i(m)}{(1-2^{2r_d})^{i+1}}+O(2^{-2m})\quad &{\rm{if}}\quad r_{d}\ne 0; \\
\frac{m^{i+1}}{i+1}+\pi_{i}(m)\quad &{\rm{if}}\quad r_{d}= 0,
\end{cases}
\]
where $P_{i}(m)$ and $\pi_{i}(m)$  are polynomials in $m$ of degree $i.$
This insight allows us to write
\EQ{keyS1}{
S_{1}(r_{d})=\begin{cases}
P_{d-3}^{\hat{\boldsymbol k},p,r_{d}}(m)+O\left(\frac{1}{2^{2m}}\right)\quad &{\rm{if}}\quad r_{d}\ne 0; \\
P_{d-2}^{\hat{\boldsymbol k},p,0}(m)\quad &{\rm{if}}\quad r_{d}= 0,
\end{cases}
}
where
\EQ{poly1}{
P_{d-3}^{\hat{\boldsymbol k},p,r_{d}}(m)=\sum_{i=1}^{d-3} \frac{a_{i}^{\hat{\boldsymbol k},p-r_d}}{(1-2^{2r_d})^{i+1}}P_i(m)\quad {\rm{and}}\quad P_{d-2}^{\hat{\boldsymbol k},p,0}(m)=
\sum_{i=1}^{d-3} a_{i}^{\hat{\boldsymbol k},p}\left(\frac{m^{i+1}}{i+1}+\pi_{i}(m)\right)
}
are  polynomials in $m$ of degree $d-3$ and $d-2$ respectively.
 \\
 
For the second sum $S_{2}(r_d)$  we can bring the identities  (\ref{v2}), (\ref{v3}) and  (\ref{v4}) together to give
\[
\sum_{\ell_d=1}^{m-d}\frac{(m-\ell_d)^{d-2}}{ 2^{2\ell_d (r_d-1)}}=\begin{cases}
\frac{P_{d-2}(m)}{(1-2^{2(r_{d}-1)})^{d-1}}+O(2^{-2m})\quad &{\rm{if}}\quad r_{d}> 1; \\
\frac{m^{d-1}}{d-1}+\pi_{d-2}(m)\quad &{\rm{if}}\quad r_{d}= 1;\\
C_{d,d-2}2^{2m}+\left(\frac{4}{3}\right)^{d-1}P_{d-2}^{*}(m)\quad &{\rm{if}}\quad r_{d}= 0,
\end{cases}
\]
and so deduce that
\[
S_{2}(r_d)=\delta_{p-r_d,d-1} (d-1)\left(k_{1}^2\cdots k_{d-1}^2\right)\begin{cases}
\frac{P_{d-2}(m)}{2^{2m}(1-2^{2(r_{d}-1)})^{d-1}}+O(2^{-4m})\quad &{\rm{if}}\quad r_{d}> 1; \\
\frac{m^{d-1}}{2^{2m}(d-1)}+\frac{\pi_{d-2}(m)}{2^{2m}}\quad &{\rm{if}}\quad r_{d}= 1;\\
C_{d,d-2}+\left(\frac{4}{3}\right)^{d-1}\frac{P_{d-2}^{*}(m)}{2^{2m}}\quad &{\rm{if}}\quad r_{d}= 0.
\end{cases}
\]
 
  Isolating the dominant  term from those exhibiting faster decay we have that
 \EQ{keyS2}{
S_{2}(r_d)=\delta_{p-r_d,d-1} \left(k_{1}^2\cdots k_{d-1}^2\right)\begin{cases}
\frac{C_{d} m^{d-2}}{2^{2m}}+O\left(\frac{m^{d-3}}{2^{2m}}\right)\quad &{\rm{if}}\quad r_{d}> 1; \\
\frac{m^{d-1}}{2^{2m}}+O\left(\frac{m^{d-2}}{2^{2m}}\right)\quad &{\rm{if}}\quad r_{d}= 1;\\
C_{d}^{\prime}+O\left(\frac{m^{d-2}}{2^{2m}}\right)\quad &{\rm{if}}\quad r_{d}= 0,
\end{cases}
}
where  $C_{d}^{\prime}= (d-1)C_{d,d-2}.$ For the third sum $S_{3}(r_d)$  we can use the same approach as above, with $d-3$ replacing $d-2,$ to deduce that
\EQ{keyS3}{
S_{3}(r_d)=C(\hat{\boldsymbol k},p-r_d,d-1)\begin{cases}
\frac{C_{d-1} m^{d-3}}{2^{2m}}+O\left(\frac{m^{d-4}}{2^{2m}}\right)\quad &{\rm{if}}\quad r_{d}>1; \\
\frac{m^{d-2}}{(d-2)2^{2m}}+O\left(\frac{m^{d-3}}{2^{2m}}\right)\quad &{\rm{if}}\quad r_{d}=1;\\
C_{d-1}^{\prime}+O\left(\frac{m^{d-3}}{2^{2m}}\right)\quad &{\rm{if}}\quad r_{d}=0.
\end{cases}
}

We now bring our findings (\ref{keyS1}),(\ref{keyS2}),(\ref{keyS3}) together, where again we isolate the dominant terms from those with faster decay to provide
\[
\begin{aligned}
\sigma_{d,p}(m,\boldsymbol k)&=\sum_{r_d=0}^{p}\frac{p! k_{d}^{2r_d}}{(p-r_{d})! r_{d}!}\left(S_{1}(r_{d})+S_{2}(r_{d})+S_{3}(r_{d})\right)+O\left(\frac{m^{d-3}}{2^{2m}}\right)\\
&=P_{d-2}^{\hat{\boldsymbol k},p,0}(m)+\sum_{r_d=1}^{p}\frac{p! k_{d}^{2r_d}}{(p-r_{d})! r_{d}!}P_{d-3}^{\hat{\boldsymbol k},p,r_{d}}(m)\\
&+\delta_{p,d-1} \left(k_{1}^2\cdots k_{d-1}^2\right)
C_{d}^{\prime} +\delta_{p-1,d-1} \left(k_{1}^2\cdots k_{d-1}^2\right) \frac{k_{d}^2 p!}{(p-1)!}\frac{m^{d-1}}{2^{2m}}+O\left(\frac{m^{d-2}}{2^{2m}}\right)\\
&+C(\hat{\boldsymbol k},p,d-1)C_{d-1}^{\prime}+C(\hat{\boldsymbol k},p-1,d-1)\frac{pk_{d}^2}{d-2}\frac{m^{d-2}}{2^{2m}}+O\left(\frac{m^{d-3}}{2^{2m}}\right).
\end{aligned}
\]

By inspection we observe that the above can be expressed as
\EQ{restdone}{
\sigma_{d,p}(m,\textbf{k})=\pi_{d-2}^{\textbf{k},p}(m)+\delta_{p,d} \frac{dk_{1}^{2}\cdots k_{d}^{2}m^{d-1}}{2^{2m}}+
 \frac{C(\textbf{k},p,d)m^{d-2}}{2^{2m}}+O\left(\frac{m^{d-3}}{2^{2m}}\right),
}
where $\pi_{d-2}^{\textbf{k},p}(m)$ is a polynomial in $m$ of degree $d-2$ whose coefficients depend upon $\textbf{k}$ and $p.$ This completes the proof of Theorem \ref{main1}.

\subsection*{Acknowledgments}

The work of Janin Jäger was funded by the Deutsche Forschungsgemeinschaft (DFG - German research foundation) - Projektnummer: 461449252.

\end{document}